\documentclass[10pt,conference,letterpaper]{IEEEtran}
\IEEEoverridecommandlockouts

\usepackage{cite}
\usepackage{amsmath,amssymb,amsfonts}
\usepackage{bm}
\usepackage{algorithmic}
\usepackage[dvipdfmx]{graphicx}
\usepackage{textcomp}
\usepackage{xcolor}
\def\BibTeX{{\rm B\kern-.05em{\sc i\kern-.025em b}\kern-.08em
    T\kern-.1667em\lower.7ex\hbox{E}\kern-.125emX}}
\begin{document}

\title{Annealing-based approach to solving partial differential equations
\thanks{This work was supported by a JSPS KAKENHI Grant Numbers JP23H04499 and JP25H01522.}
}

\author{\IEEEauthorblockN{Kazue Kudo}
\IEEEauthorblockA{\textit{Department of Computer Science} \\
\textit{Ochanomizu University}\\
Tokyo, Japan\\
https://orcid.org/0000-0002-0060-7772}
}
\maketitle

\begin{abstract}
Solving partial differential equations (PDEs) using an annealing-based approach involves solving generalized eigenvalue problems. 
Discretizing a PDE yields a system of linear equations (SLE). 
Solving an SLE can be formulated as a general eigenvalue problem, which can be transformed into an optimization problem with an objective function given by a generalized Rayleigh quotient.
The proposed algorithm requires iterative computations. 
However, it enables efficient annealing-based computation of eigenvectors to arbitrary precision without increasing the number of variables. 
Investigations using simulated annealing demonstrate how the number of iterations scales with system size and annealing time.
Computational performance depends on system size, annealing time, and problem characteristics.
\end{abstract}


\section{Introduction}
\label{sec:intro}

Partial differential equations (PDEs) are ubiquitous in science and engineering, with applications spanning numerous fields. 
While analytical solutions to PDEs are often hard to find, computational approaches are indispensable for obtaining an approximate solution. 
PDEs arising from practical problems often involve large-scale systems and numerous variables, demanding substantial computational resources. 
To address these computational challenges, several quantum algorithms offering exponential speedups have been proposed~\cite{Cao2013,Childs2021}. 
These algorithms rely on fault-tolerant quantum computers, which are not yet technologically feasible.
Using those algorithms is currently an unrealistic choice. 
As a more practical alternative, variational quantum algorithms (VQAs) have emerged as a promising approach for near-term quantum devices. 
Several VQAs have been proposed for solving PDEs~\cite{Liu2021,Sato2021,Liu2022,Demirdjian2022}.

Discretizing a PDE leads to a system of linear equations (SLE).
Several VQAs have been proposed for solving SLE~\cite{Liu2021,Sato2021,Liu2022,Demirdjian2022,Bravo2023,Sato2023}.
The problem of solving an SLE can be formulated as a generalized eigenvalue problem.
VQAs for generalized eigenvalue problems involve optimizing a generalized Rayleigh quotient. 
Solving optimization problems is also the field where quantum annealing has an advantage.
Several annealing-based algorithms have been proposed for solving eigenvalue problems~\cite{Teplukhin2019,Teplukhin2020,Krakoff2022}.

This paper presents a method for solving PDEs using an Ising machine.
Ising machines are special-purpose computers designed to efficiently solve combinatorial optimization problems.
They have been realized using various technologies, including a quantum annealer~\cite{Johnson2011}, digital processors based on simulated annealing (SA)~\cite{Yamaoka2016,Kihara2017,Aramon2019,Okuyama2019,Yamamoto2020}, those based on simulated bifurcation~\cite{Goto2019,Tatsumura2020,Goto2021,Hidaka2023}, and coherent Ising machines~\cite{Wang2013,Marandi2014,Inagaki2016,McMahon2016,Hamerly2019}.
Although initially designed to solve combinatorial optimization problems, they have applications in various fields.
This work demonstrates their potential for solving PDEs and expands their application fields.

For solving a problem using an Ising machine, the problem is represented by the quadratic unconstrained binary optimization (QUBO) formulation.
In the QUBO formulation, the solution to a PDE is expressed as a linear combination of binary variables.
The number of binary variables, which depends on the desired precision and the discretization mesh size, is a critical factor in QUBO-based methods.
This work employs an iterative method, inspired by Ref.~\cite{Krakoff2022}, that uses a relatively small number of variables.
In a usual QUBO formulation, increasing precision typically requires more binary variables. 
However, the proposed algorithm enables the computation of eigenvectors to arbitrary precision without increasing the number of variables.
The key idea is to improve accuracy by refining the mesh-size scale while keeping the number of divisions fixed.
Another advantage of the algorithm is that the desired accuracy becomes available even if the performance of the annealing procedure is imperfect.
This paper introduces a modified algorithm that improves upon the original one and proposes a more efficient precision update scheme.

The proposed algorithm, when implemented on Ising machines, has the potential to be a promising PDE solver, although the algorithm's heuristic nature may limit its accuracy. 
Numerical simulations using simulated annealing (SA) demonstrate how the number of iterations scales with the number of variables.
The scaling behavior varies across different PDEs.  
Numerical simulations also reveal how the success rate of obtaining accurate solutions depend on the number of variables.
These insights can inform the selection of optimal parameters for solving a PDE using an Ising machine. 

\section{Methods}
\label{sec:method}

\subsection{Generalized eigenvalue problem}
We consider the SLE arising from the discretization of a PDE with Dirichlet boundary conditions:
\begin{align}
 K\bm{u} = \bm{f},
\label{eq:SLE}
\end{align}
where $K$ is a symmetric positive-definite matrix, and $\bm{u}$ and $\bm{f}$ are real vectors.
While SLEs can be defined in complex spaces, we focus on real-valued systems for simplicity.
The goal is to determine the solution $\bm{u}$ for given $K$ and $\bm{f}$.

The SLE can be transformed into the form of a generalized eigenvalue problem, which is generally written as
\begin{align}
 A\bm{v} = \lambda B\bm{v},
\label{eq:GEP}
\end{align}
where $A$ and $B$ are symmetric matrices, and $B$ is positive definite.
Here, $\lambda$ and $\bm{v}$ denote the generalized eigenvalue and eigenvector, respectively.
In this case, the smallest generalized eigenvalue $\lambda_{\rm min}$ is equivalent to the minimum generalized Rayleigh quotient:
\begin{align}
 \lambda_{\rm min} = \min_{\bm w}
\frac{\bm{w}^\top A\bm{w}}{\bm{w}^\top B\bm{w}},
\end{align}
where $\bm{w}$ is a unit vector.
The eigenvector $\bm{v}_{\rm min}$ corresponding to $\lambda_{\rm min}$ minimizes the generalized Rayleigh quotient.

Here, we rewrite \eqref{eq:SLE} as
\begin{align}
 \tilde{K}\bm{u} = \tilde{\bm f},
\label{eq:SLE.1}
\end{align}
where $\tilde{K}=K/\|\bm{f}\|$, and $\tilde{\bm f}=\bm{f}/\|\bm{f}\|$ is a unit vector.
Equation~\eqref{eq:SLE.1} can be formulated as
\begin{align}
 -\tilde{\bm f}\tilde{\bm f}^\top\bm{v}=\lambda\tilde{K}\bm{v},
\label{eq:GEP.1}
\end{align}
which has the same form as \eqref{eq:GEP}.
Comparing \eqref{eq:GEP} and \eqref{eq:GEP.1}, we identify $A=-\tilde{\bm f}\tilde{\bm f}^\top$ and $B=\tilde{K}$.
Since $-\tilde{\bm f}\tilde{\bm f}^\top$ is a rank-$1$ matrix, equation \eqref{eq:GEP.1} has only one nonzero eigenvalue, which is negative and corresponds to the smallest eigenvalue $\lambda_{\rm min}$, and the other eigenvalues are zero.
Dividing both sides of \eqref{eq:GEP.1} by $\tilde{\bm f}^\top\bm{v}_{\rm min}$, which is nonzero as $\lambda_{\rm min}\neq 0$, we obtain
\begin{align}
 \tilde{K}\left(-\frac{\lambda_{\rm min}\bm{v}_{\rm min}}{\tilde{\bm f}^\top\bm{v}_{\rm min}}\right)=\tilde{\bm f}.
\label{eq:SLE.2}
\end{align}
Combining \eqref{eq:SLE.1} and \eqref{eq:SLE.2} results in
\begin{align}
 \tilde{K}\left(\bm{u}+\frac{\lambda_{\rm min}\bm{v}_{\rm min}}{\tilde{\bm f}^\top\bm{v}_{\rm min}}\right)=\bm{0}.
\end{align}
Considering that $\tilde{K}$ is positive definite, i.e., invertible, we arrive at the solution to \eqref{eq:GEP.1}:
\begin{align}
 \bm{u} = -\frac{\lambda_{\rm min}\bm{v}_{\rm min}}{\tilde{\bm f}^\top\bm{v}_{\rm min}}.
\label{eq:u}
\end{align}

\subsection{Algorithm}

Here, we review the core of the algorithm proposed in Ref.~\cite{Krakoff2022}, aiming to solve the generalized eigenvalue problem:
\begin{align}
 A\bm{v}=\lambda B\bm{v},
\end{align}
where $A,B\in\mathbb{R}^{n\times n}$ are symmetric matrices, and $B\in\mathbb{R}^{n\times n}$ is positive definite.
The eigenvector $\bm{v}\in\mathbb{R}^n$ is approximated as $\bm{v}=(I_n\otimes\bm{p})\bm{q}$, where $\bm{p}=\left(-1,\frac12,\frac1{2^2},\ldots, \frac1{2^{b-1}}\right)$ is a row vector, and $\bm{q}=\{0,1\}^{nb}$ is a binary vector.
The parameter $b$ determines the precision of the initial guess.
The explicit form of $\bm{v}$ is written as
\begin{align}
 \bm{v}=\left[\begin{array}{c|c|c|c}
-1,\frac12,\cdots,\frac{1}{2^{b-1}} & 0 & \cdots & 0
\\ \hline
0 & \ddots &\ddots  & \vdots
\\ \hline
\vdots &\ddots  & \ddots & 0	 
\\ \hline
0 & \cdots & 0 & -1,\frac12,\cdots,\frac{1}{2^{b-1}}
\end{array}\right]\bm{q}.
\label{eq:var}
\end{align}
Each element of $\bm{v}$ takes a discrete value in the range $[-1,\; 1-\frac1{2^{b-1}}]$.

This algorithm consists of two stages.
In the first stage, which we call the initial guess stage, an initial guess for the solution is obtained with a precision of $\frac1{2^{b-1}}$.
In the second stage, which we call the iterative descent stage, an iterative optimization procedure is performed to refine the solution with increasing precision.

\subsubsection{Initial guess stage}

In the initial guess stage, the following two steps are repeated until $\lambda$ converges:
\begin{itemize}
 \item[(1)] Determine a tentative eigenvector $\bm{v}^*$ by minimizing the generalized Rayleigh quotient:
\begin{align}
 \bm{v}^* = \underset{\bm v}{\operatorname{argmin}} \left[
\bm{v}^\top (A-\lambda B)\bm{v}\right].
\label{eq:ini1}
\end{align}
 \item[(2)] Update the eigenvalue estimate:
\begin{align}
  \lambda=\frac{\bm{v}^{*\top} A\bm{v}^*}{\bm{v}^{*\top} B\bm{v}^*}.
\label{eq:ini2}
\end{align}
\end{itemize}
The initial value of $\lambda$ is set to $\lambda=\bm{x}^\top A\bm{x}$, where $\bm{x}$ is a random unit vector.
To solve the minimization problem in Step 1, an annealing-based approach (i.e., an Ising machine or SA) is used to find $\bm{q}^*$ that minimized the QUBO:
\begin{align}
 \bm{q}^* = \underset{\bm q}{\operatorname{argmin}} \left[
\bm{q}^\top (A-\lambda B)\bm{q}\right].
\end{align}
The unnormalized eigenvector $\bm{v}=(I_n\otimes\bm{p})\bm{q}^*$ is then calculated, and normalized to obtain the tentative eigenvector $\bm{v}^*=\bm{v}/\|\bm{v}\|$.
The convergence of this inexact parametric algorithm for mixed-integer fractional programs has been proved for a good approximation of the solution~\cite{Zhong2014}.
While Ref.~\cite{Krakoff2022} proposed several additional techniques for obtaining better solutions, we focus on this straightforward approach in this study.

\subsubsection{Iterative descent stage}

The formal objective function to be minimized in the iterative descent stage is $g(\bm{x})=\bm{x}^\top(A-\lambda B)\bm{x}$.
Its Taylor expansion around $\bm{x}_0$ is written as
\begin{align}
  g(\bm{x})=g(\bm{x}_0)+
\bm{d}^\top\bm{\nabla}g(\bm{x}_0)
+\frac12\bm{d}^\top(\nabla^2 g(\bm{x}_0))\bm{d}
+ O(\|\bm{d}\|^3),
\label{eq:f}
\end{align}
where $\bm{d}=\bm{x}-\bm{x}_0$.
The gradient $\bm{\nabla}g$ and Hessian $\nabla^2 g$ are given by
\begin{align}
 \bm{\nabla}g &=2(A-\lambda B)\bm{x}, \\
 \nabla^2 g &=2(A-\lambda B).
\end{align}
The core of the iterative descent stage involves computing the update direction $\bm{d}$ that minimizes $g$ and then updating the solution: $\bm{x}_0\leftarrow\bm{x}_0 +t\bm{d}$, where $t$ is real. 

To obtain the eigenvector corresponding to a given eigenvalue $\lambda$, we seek a difference vector $\bm{d}^*$ that minimizes the objective function:
\begin{align}
 \bm{d}^*=  \underset{\bm d}{\operatorname{argmin}}\left[
2\bm{v}^\top(A-\lambda B)\bm{d}+\bm{d}^\top(A-\lambda B)\bm{d}\right],
\end{align}
where $\bm{v}$ is a tentative eigenvector.
Since the eigenvector is supposed to be a unit vector, vector $\bm{v}$ should move on the unit sphere.
Therefore, the difference should be orthogonal to $\bm{v}$.
Orthogonalizing $\bm{d}^*$, we obtain the adjusted difference vector $\bm{d}=\bm{d}^*-(\bm{v}^\top\bm{d}^*)\bm{v}$.
The eigenvector is then updated as $\bm{v}\leftarrow\bm{v} +t\bm{d}$, where the parameter $t$ is determined to minimize
 \begin{align}
 g(\bm{v}+t\bm{d})=(\bm{v}+t\bm{d})^\top(A-\lambda B)(\bm{v}+t\bm{d}).
\label{eq:ft}
\end{align}
This quadratic function of $t$ has a minimum only if $\bm{d}^\top(A-\lambda B)\bm{d}>0$.
Otherwise, we set $t=1$.
To ensure that $t\ge 1$, we define the parameter as
\begin{align}
 t=
\begin{cases}
\mathrm{sgn}(b) \max\left( \frac{|b|}{a},\; 1\right) & a>0, \\
1 & a\le 0,
\end{cases}
\label{eq:t}
\end{align}
where $a=\bm{d}^\top(A-\lambda B)\bm{d}$ and $b=-\bm{v}^\top(A-\lambda B)\bm{d}$.
If $t\ge 1$, $\bm{v} +t\bm{d}$ may overshoot the optimal solution, indicating that the discretization mesh size is not small enough to produce a better solution.
In this case, we refine the mesh size instead of updating $\bm{v}$.
When the mesh size reaches $\epsilon_0/2^{b-1}$, where $\epsilon_0$ is a tolerance precision parameter, the algorithm terminates.

In principle, the above procedure can solve the generalized eigenvalue problem to arbitrary precision.
The original algorithm proposed in Ref.~\cite{Krakoff2022} often provides solutions better than single precision but worse than double precision.
This issue stems from the mesh size downscaling condition.
Even without overshooting, the annealing procedure's performance can limit solution improvement, leading to unnecessary mesh refinement.
To avoid such a situation, the annealing procedure should be iterated several times at a given mesh size when the solution shows no improvement.

The procedure in the iterative descent stage is summarized as follows.
First, set the initial mesh size: $r=r_{\rm ini}$.
Then, the following procedure is repeated until $r$ reaches to the tolerance precision $\epsilon_0$ or the eigenvalue $\lambda$ converges.
\begin{itemize}
 \item[(1)] Using an annealing-based approach, solve the QUBO problem to obtain the binary vector $\bm{q}^*$:
\begin{align}
 \bm{q}^* = \underset{\bm q}{\operatorname{argmin}} \left[
2r\bm{v}^\top C(I_n\otimes\bm{p})\bm{q}
+r^2\bm{q}^\top (C\otimes\bm{p}^\top\bm{p})\bm{q}\right],
\end{align}
where $C=A-\lambda B$.
 \item[(2)] Calculate the adjusted difference vector 
       $\bm{d}=\bm{d}^*-(\bm{v}^\top\bm{d}^*)\bm{v}$, where $\bm{d}^*=(I_n\otimes\bm{p})r\bm{q}^*$.
 \item[(3)] Calculate the tentative eigenvector as $\bm{v}^* = \bm{v} +t\bm{d}$, where $t$ is defined by \eqref{eq:t}.
 \item[(4)] Update the tentative eigenvalue estimate:
\begin{align}
  \lambda^*=\frac{\bm{v}^{*\top} A\bm{v}^*}{\bm{v}^{*\top} B\bm{v}^*}.
\end{align}
 \item[(5a)] If $\lambda^*<\lambda$, accept the updates, i.e., 
       $\lambda=\lambda^*$ and $\bm{v}=\bm{v}^*$, 
       and return to Step 1.
 \item[(5b)] If $\lambda^*\ge \lambda$ and the number of consecutive rejections is less than a threshold $n_{\rm rpt}$, accept only the eigenvector update, i.e., $\bm{v}=\bm{v}^*$, and return to Step 1.
 \item[(5c)] Otherwise, reduce the mesh size: $r\leftarrow\eta r$, where $\eta$ is real, and return to Step 1.
\end{itemize}
The precision update rate $\eta$ was a constant in the original algorithm proposed in Ref.~\cite{Krakoff2022}.
In this study, the update rule is modified so that the updated precision is comparable to mesh size, i.e., $\eta=2^{-b+1}$.

\section{Results}
\label{sec:result}

This study investigates the computational efficiency of the proposed method for solving PDEs, focusing on the number of iterations required to achieve a given precision. 
Here, SA was employed as the underlying annealing procedure.
While using an Ising machine can potentially accelerate computations, it does not directly impact the achievable precision. 

To assess the success rate of obtaining an accurate solution, we consider problems with known exact solutions.
The success rate was defined as the proportion of trials where the obtained solution's root mean squared error (RMSE) is smaller than a specified threshold $\epsilon_1$.
Here, the true solution, used as a reference for RMSE calculation, was determined by the solution of the corresponding SLE and was obtained by a classical linear equation solver. 

This study investigates three types of Poisson equations of the form:
\begin{align}
 -\nabla^2u(\bm{x})=f(\bm{x}).
\label{eq:Poi}
\end{align}
The first and second ones are one-dimensional symmetric and asymmetric Poisson equations with $f(\bm{x})$ given by 
\begin{align}
 f_{\rm sym}(x)&=-12x^2+12x-2,
\label{eq:sym}\\
 f_{\rm asym}(x)&=-12x^2+18x-5.5,
\label{eq:asym}
\end{align}
respectively, where $x\in [0,1]$.
The boundary conditions are $u(x)=0$ at $x=0,1$.
The third one is the two-dimensional Poisson equation with $f(\bm{x})$ given by
\begin{align}
  f_{\rm 2D}(x,y)=-2x(x-1)-2y(y-1),
\label{eq:2D}
\end{align}
where $x,y\in [0,1]$.
The boundary conditions are $u(x,y)=0$ at $x=0,1$ or $y=0,1$.
The exact solutions for these equations are $u(x)=x^2(x-1)^2$ for the one-dimensional symmetric case, $u(x)=x(x-1)(2x-1)(2x-3)/4$ for the one-dimensional asymmetric case, and $u(x,y)=xy(x-1)(y-1)$ for the two-dimensional case.

The system size of the SLE corresponding to the Poisson equation is determined by the dimension of vector $\bm{u}$ in \eqref{eq:SLE.1}.
When the range $[0,1]$ is discretized into $n+1$ points in each dimension, the system size is $n$ and $n^2$ for one- and two-dimensional problems, respectively.

Unless otherwise specified, the following parameters were used: $r_{\rm ini}=\eta=2^{-b+1}$, $\epsilon_0=\epsilon_1=10^{-8}$, and $n_{\rm rpt}=0$.
For each experimental condition, $10^3$ trials were performed, and results were averaged over the executions.

In this study, we performed SA using the open-source Python library OpenJij~\cite{OpenJij}. 
The SA sampler in this library automatically determines a geometric temperature schedule.

\subsection{Dependence on precision}
\label{sec:b_dependence}

\begin{figure}[htbp]
\centering
 \includegraphics[width=0.25\textwidth]{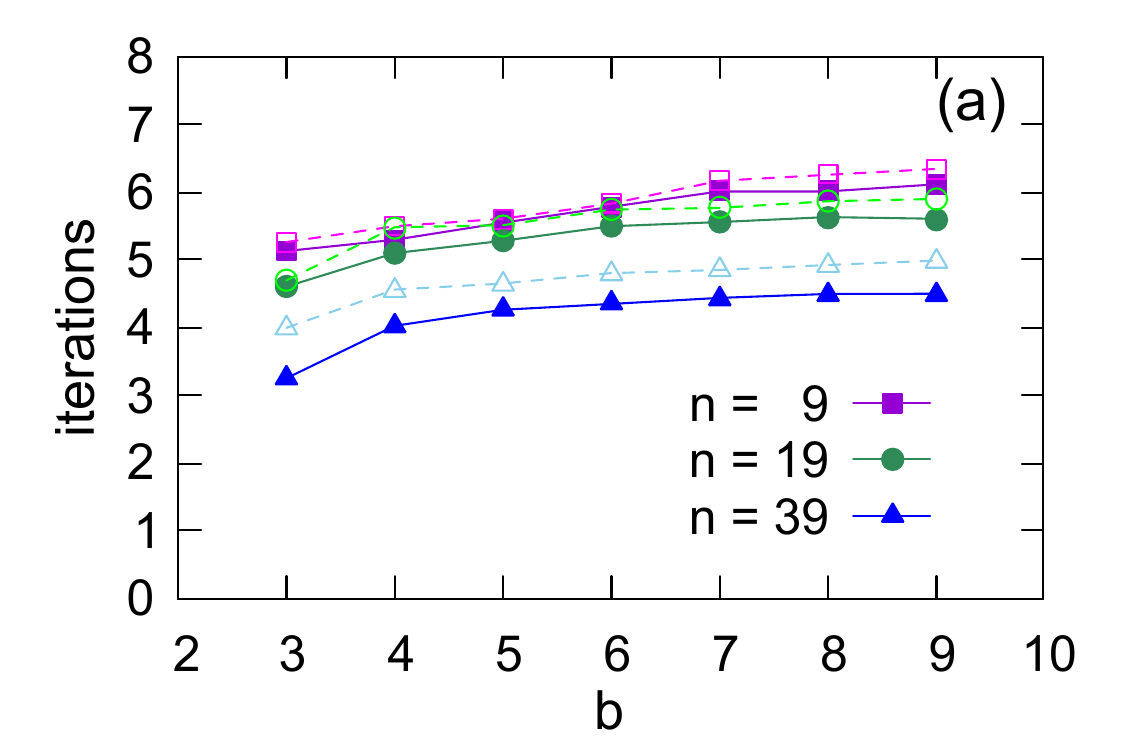}\\
 \includegraphics[width=0.25\textwidth]{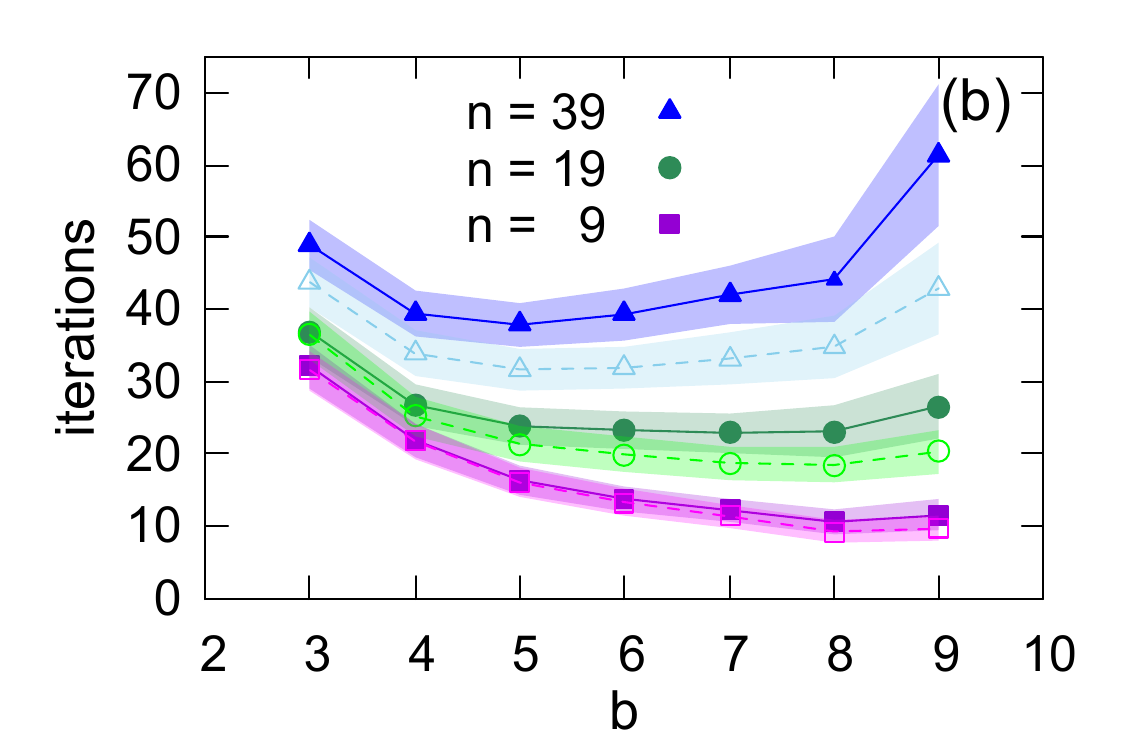}%
 \includegraphics[width=0.25\textwidth]{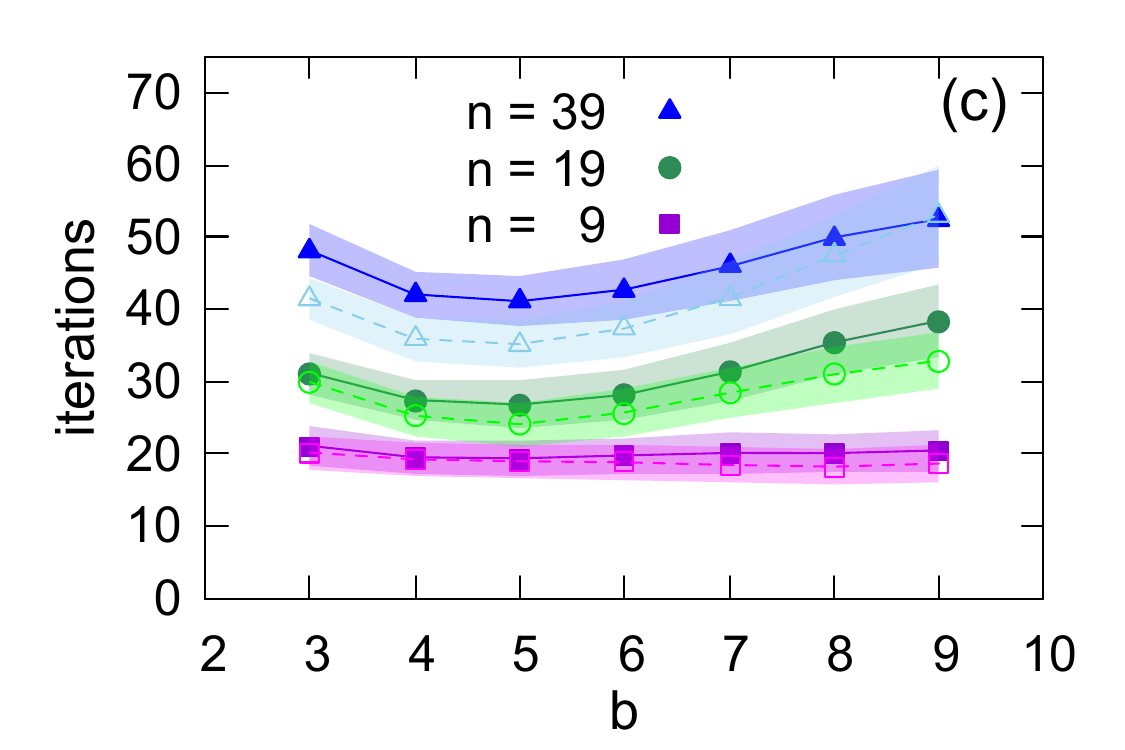}%
\caption{\label{fig:tot-itr}
The number of iterations required for solving the one-dimensional symmetric Poisson equation \eqref{eq:sym} as a function of the precision parameter $b$.
(a) the initial guess stage. (b) iterative descent stage with $\eta=2^{-b+1}$. (c) iterative descent stage with $\eta=0.1$.
The dimensions of solution vectors are $n=9$, $n=19$, and $n=39$.
Filled symbols with solid lines correspond to short annealing times, and open symbols with dashed lines correspond to long annealing times.
Symbols and shaded ranges represent the average and standard deviation, respectively.
Symbol size is proportional to the success rate.}%
\end{figure}

Figure~\ref{fig:tot-itr} illustrates the number of iterations required to solve the one-dimensional symmetric Poisson equation as a function of the precision parameter $b$ at the initial guess and iterative descent stages.
Filled symbols with solid lines correspond to $N_{\rm step}=10^3$, while open symbols with dashed lines correspond to $N_{\rm step}=10^4$, where $N_{\rm step}$ is the number of temperature steps in the annealing process and corresponds to annealing time.

At the initial guess stage, as shown in Fig.~\ref{fig:tot-itr}(a), the number of iterations is relatively insensitive to the precision parameter $b$.
The standard deviation, which is not depicted, typically ranges from 0.5 to 1.0.
The number of iterations is smaller for shorter annealing times and decreases as the system size increases.
These trends are different from those in Figs.~\ref{fig:tot-itr}(b) or (c), which reflects the difference in objective functions.

In contrast, the iterative descent stage exhibits a stronger dependence on $b$. 
Figure~\ref{fig:tot-itr}(b) shows the results for a $b$-dependent mesh reduction rate ($\eta=2^{-b+1}$), while Fig.~\ref{fig:tot-itr}(c) shows the results for a constant mesh reduction rate ($\eta=0.1$). 
The $b$-dependent mesh reduction generally leads to fewer iterations, especially for smaller $n$.
However, for large $n$, the number of iterations increases with $b$ rapidly.

In Figs.~\ref{fig:tot-itr}(b) and \ref{fig:tot-itr}(c), the symbol size is proportional to the success rate. 
For example, the smaller symbol size for $n=39$ and $b=8$ indicates a low success rate.
In both the $b$-dependent and constant mesh reduction cases, larger system sizes require more iterations.
Although increasing the annealing time reduces the number of iterations, the improvement is limited.
Even with a tenfold increase in annealing time, the number of iterations decreases by less than 20\%. 

\begin{figure}[htbp]
\centering
 \includegraphics[width=0.25\textwidth]{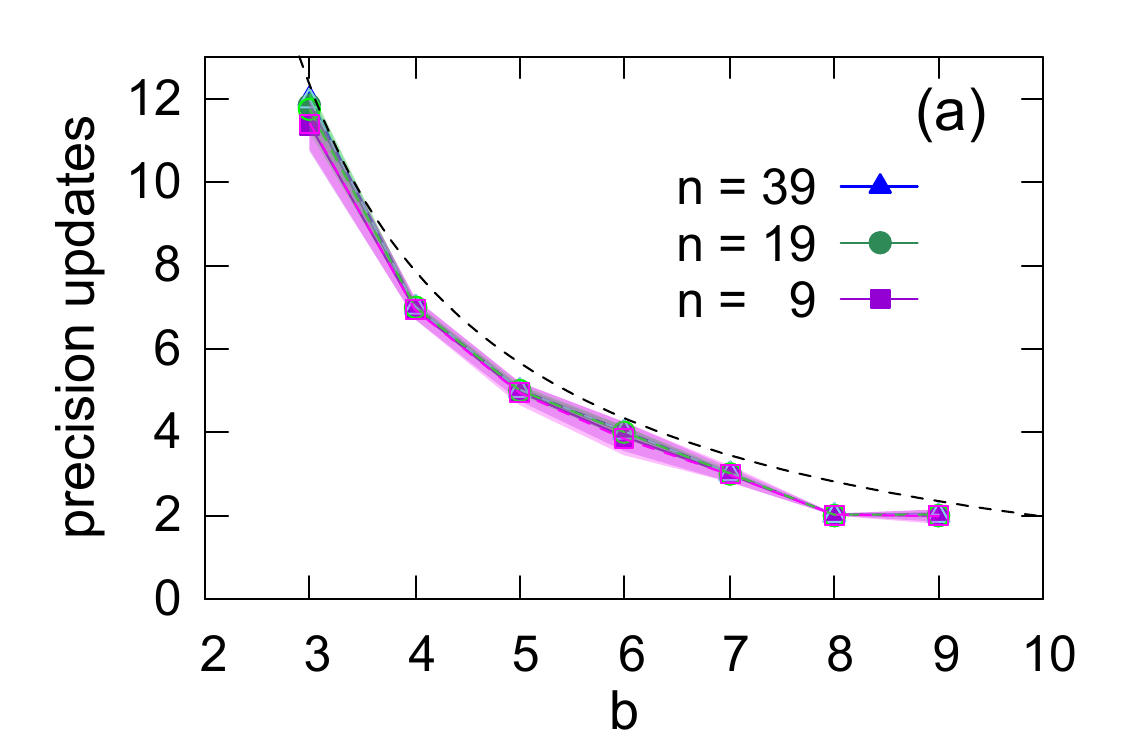}%
 \includegraphics[width=0.25\textwidth]{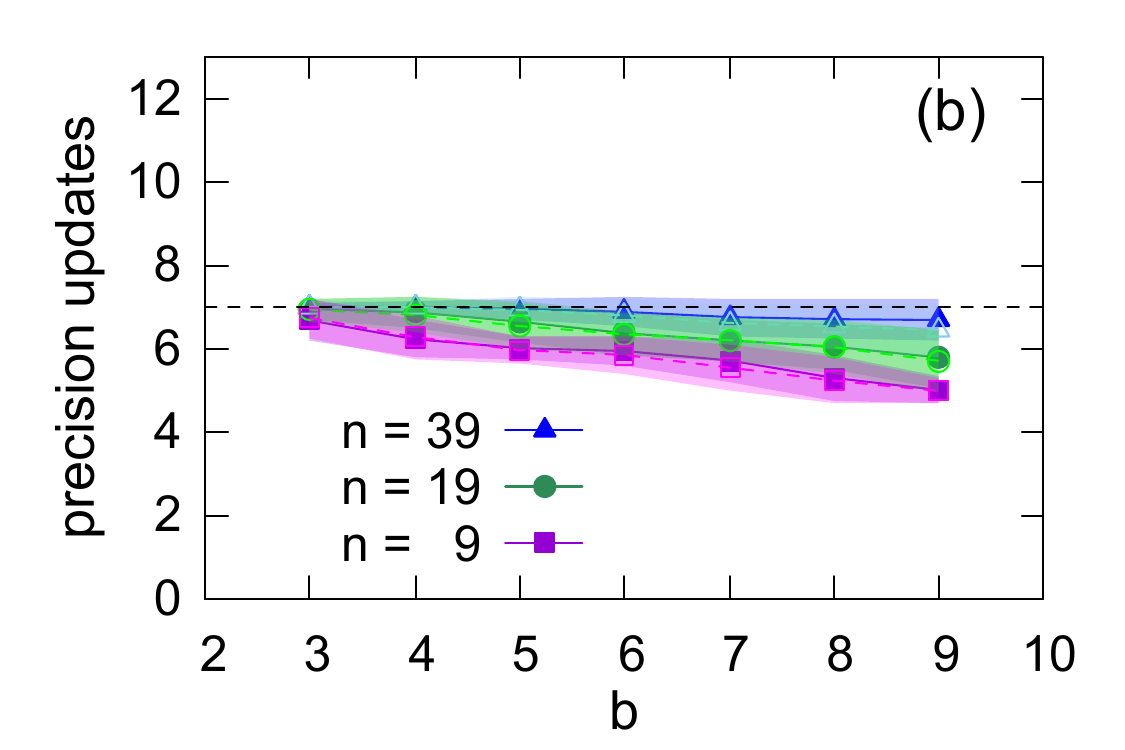}%
\caption{\label{fig:upd}
The number of precision updates required for solving the one-dimensional symmetric Poisson equation \eqref{eq:sym} as a function of the precision parameter $b$. 
(a) $b$-dependent mesh reduction rate ($\eta=2^{-b+1}$). 
(b) Constant mesh reduction rate ($\eta=0.1$).
Black dashed lines indicate the theoretical estimate from \eqref{eq:n_upd.0} for $r=\epsilon_0=10^{-8}$.
The dimensions of solution vectors are $n=9$, $n=19$, and $n=39$.
Filled symbols with solid lines correspond to short annealing times, and open symbols with dashed lines correspond to long annealing times.
Symbols and shaded ranges represent the average and standard deviation, respectively.}%
\end{figure}

Here, we estimate the number of precision updates in the iterative descent stage.
When the mesh scale reaches $r$ after $n_{\rm upd}$ updates, $r=r_{\rm ini}\eta^{n_{\rm upd}}$, which leads to
\begin{align}
    n_{\rm upd} = \frac{\log_2r - \log_2r_{\rm ini}}{\log_2\eta}.
\label{eq:n_upd.0}
\end{align}
For a constant $\eta$, $n_{\rm upd}$ is constant unless the eigenvalue $\lambda$ converges in the middle of the algorithm.
However, when $\eta=r_{\rm ini}=2^{-b+1}$, the expression simplifies to
\begin{align}
    n_{\rm upd} = \frac{\log_2r^{-1}}{b-1} - 1.
\label{eq:n_upd}
\end{align}
Figure~\ref{fig:upd} confirms this estimate.

Figures~\ref{fig:upd}(a) and \ref{fig:upd}(b) show the number of precision updates corresponding to Figs.~\ref{fig:tot-itr}(b) and \ref{fig:tot-itr}(c), respectively.
Almost no difference appears between different annealing times ($N_{\rm steps}=10^3$ and $10^4$).
The black dashed lines represent theoretical estimates: Eq.~\eqref{eq:n_upd} with $r=\epsilon_0$ for Fig.~\ref{fig:upd}(a) and Eq.~\eqref{eq:n_upd.0} with $r=\epsilon_0$ and $r_{\rm ini}=\eta=0.1$ for Fig.~\ref{fig:upd}(b).
These lines indicate the maximum possible number of updates.
In Fig.~\ref{fig:upd}(a), the data points closely align with theoretical estimates, regardless of $n$.
In contrast, Fig.~\ref{fig:upd}(b) shows a decreasing trend in the number of updates with increasing $b$ for smaller system sizes.
Data points below the dashed line indicate early convergence of the eigenvalue $\lambda$.

For smaller $n$ and $b$, the number of iterations in the iterative descent stage is strongly correlated with the number of precision updates.
However, for large $n$, the number of iterations increases with $b$ more rapidly in Fig.~\ref{fig:tot-itr}(b) than in Fig.~\ref{fig:tot-itr}(c).
This rapid increase can be attributed to the insufficient number of precision updates in the $b$-dependent case, as shown in Fig.~\ref{fig:upd}(a).
When the number of precision updates is small, rapid mesh size reduction causes difficulty in adjusting to a new small mesh size.
Moreover, the dimension of a search space increases exponentially with $b$ and $n$, which results in a worse convergence.

\begin{figure}[htbp]
\centering
 \includegraphics[width=0.25\textwidth]{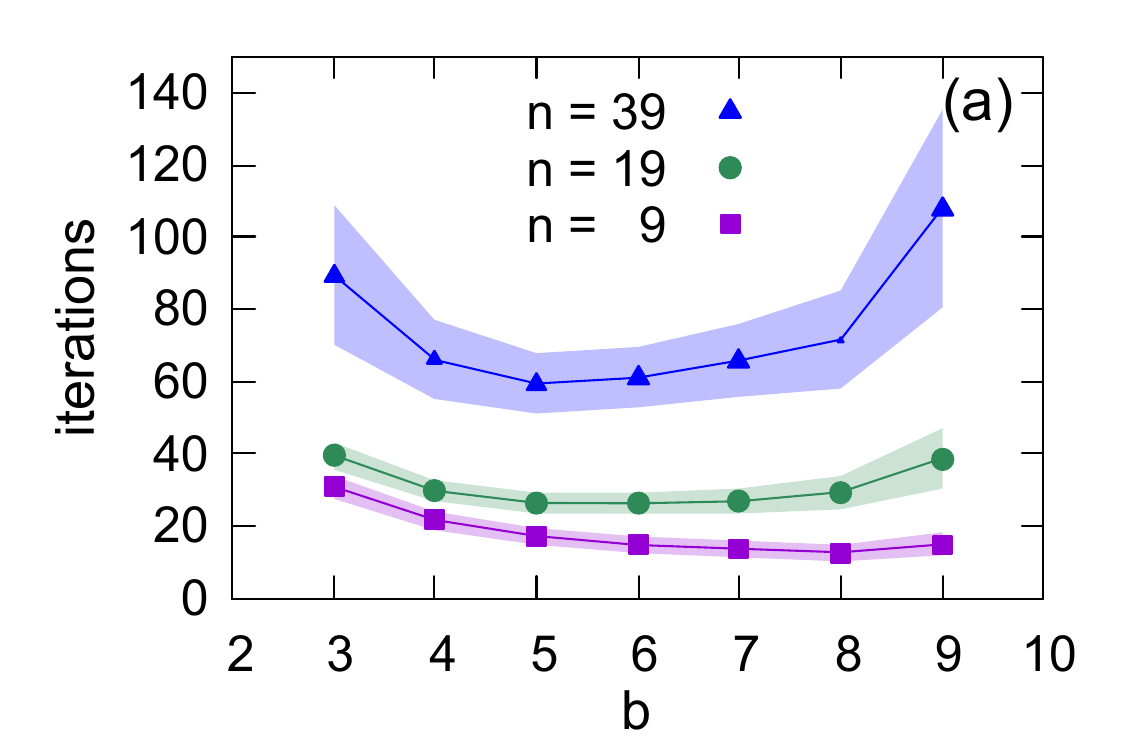}%
 \includegraphics[width=0.25\textwidth]{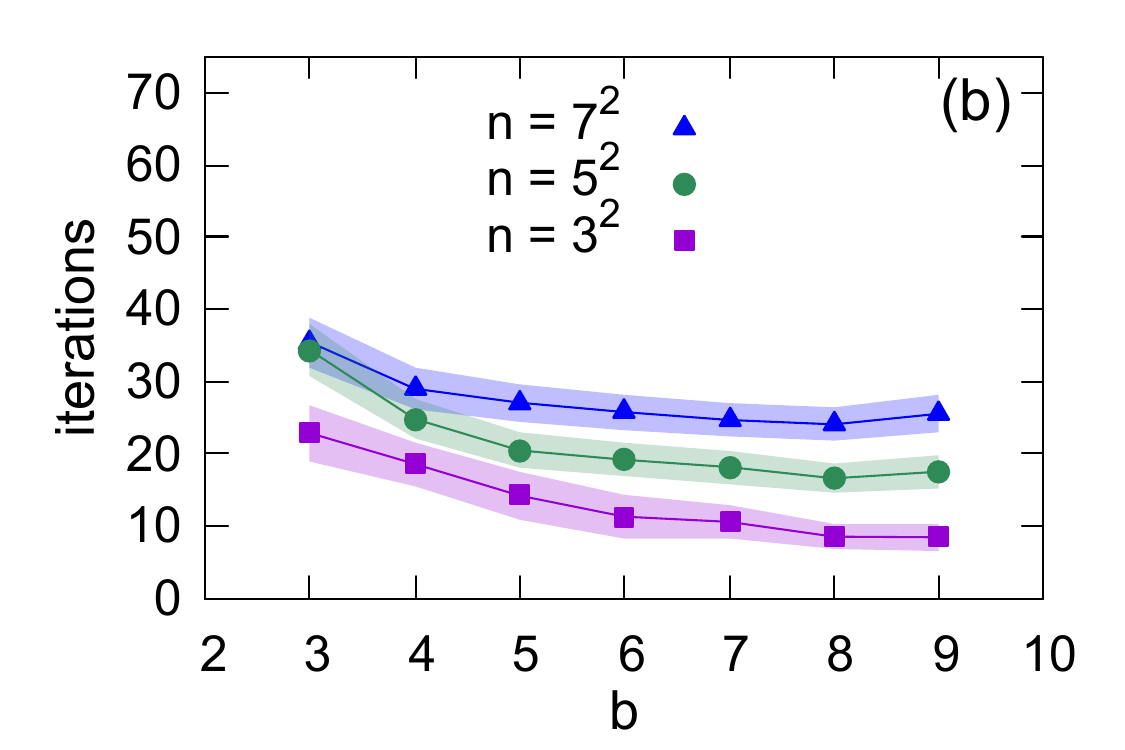}%
\caption{\label{fig:itr2}
The number of iterations at the iterative descent stage required for solving (a) the one-dimensional asymmetric Poisson equation~\eqref{eq:asym} and (b) the two-dimensional Poisson equation~\eqref{eq:2D} as a function of the precision parameter $b$.
The dimensions of solution vectors are $n=9$, $n=19$, and $n=39$ for (a) and $n=3^2$, $n=5^2$, and $n=7^2$ for (b).
Symbols and shaded ranges represent the average and standard deviation, respectively.
Symbol size is proportional to the success rate.}
\end{figure}

Figure~\ref{fig:itr2} illustrates the number of iterations at the iterative descent stage for solving (a) the one-dimensional asymmetric Poisson equation and (b) the two-dimensional Poisson equation. 
Only results for $N_{\rm step}=10^3$ are plotted.

Comparing Fig.~\ref{fig:itr2}(a) to Fig.~\ref{fig:tot-itr}(b), we observe a significantly higher number of iterations, particularly for larger $n$.
This difference can be attributed to the asymmetry of the exact solution $u(x)=x(x-1)(2x-1)(2x-3)/4$.
Figure~\ref{fig:itr2}(a) indicates that solving more complex problems generally requires more iteration.
Note that symbols for $b=4,8$ and $n=39$ are apparently small, which indicates low success rates. 

In contrast, Fig.~\ref{fig:itr2}(b) shows a much smaller number of iterations for the two-dimensional Poisson equation, even for comparable system sizes $n$ to those in Fig.~\ref{fig:tot-itr}.
This behavior can be attributed to the symmetric and convex nature of the exact solution $u(x,y)=xy(x-1)(y-1)$, which makes the optimization problem relatively easier.
It is remarkable that the number of iterations decreases with increasing $b$ even for $n=7^2$.

\subsection{Dependence on system size}
\label{sec:n_dependence}

\begin{figure}[htbp]
\centering
 \includegraphics[width=0.25\textwidth]{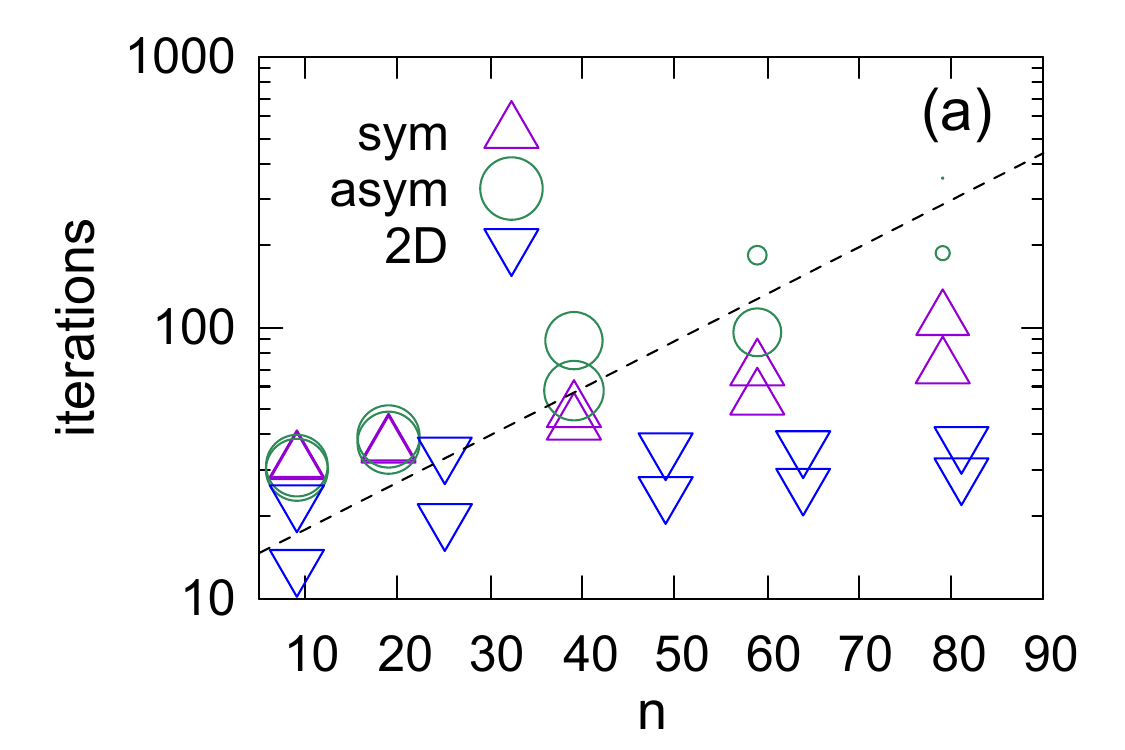}%
 \includegraphics[width=0.25\textwidth]{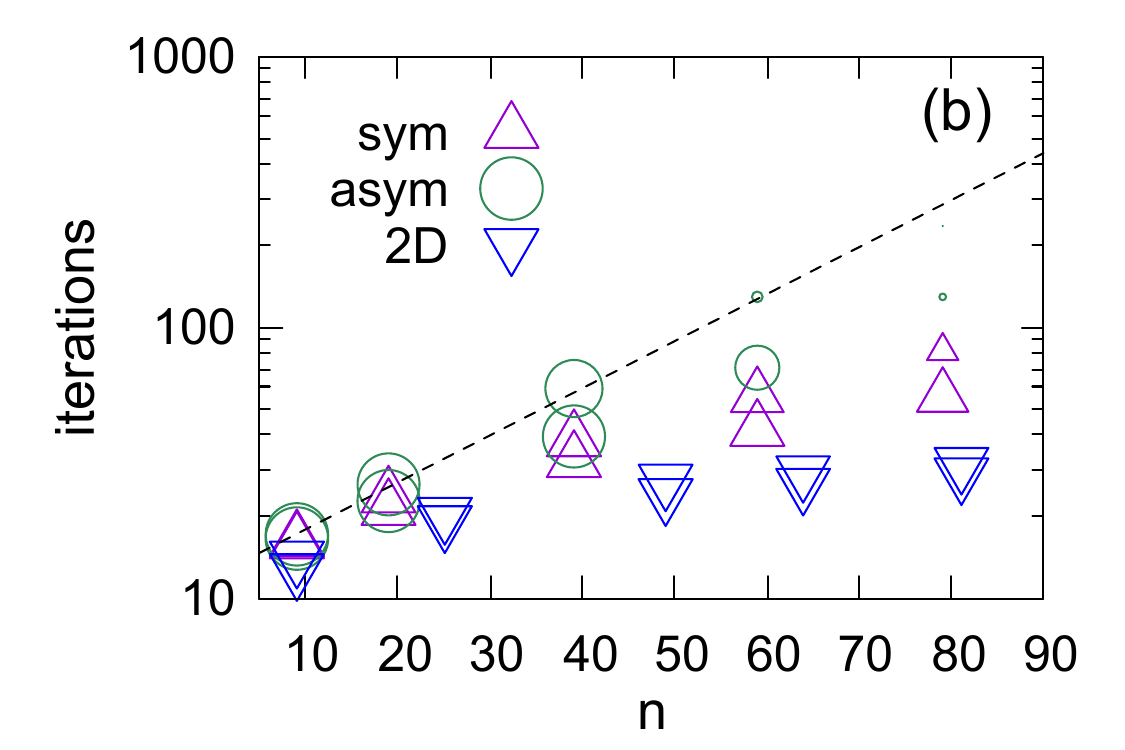}%
\caption{\label{fig:n_id}
The average number of iterations at the iterative descent stage for solving the one-dimensional symmetric, one-dimensional symmetric asymmetric, and two-dimensional Poisson equations as a function of system size $n$.
Results for $N_{\rm step}=10^3,\, 10^4$ are shown for (a) $b=3$ and (b) $b=5$.
Symbol size is proportional to the success rate.
The dashed line represents $ce^{\alpha n}$ with constants $c$ and $\alpha$.}
\end{figure}

As shown in the above, the number of iterations depends on system size and annealing time.
Figure~\ref{fig:n_id} shows the average number of iterations required at the iterative descent stage for solving the one-dimensional symmetric, one-dimensional asymmetric, and two-dimensional Poisson equations as a function of system size $n$.
The precision parameters are $b=3$ and $b=5$ in Figs.~\ref{fig:n_id}(a) and \ref{fig:n_id}(b), respectively.
Results for $N_{\rm step}=10^3$ and $10^4$ are shown in each panel, where the lower one of two points with the same symbol type at the same $n$ corresponds to $N_{\rm step}=10^4$.
Note that the symbol size is proportional to the success rate, and several points are too small to recognize, particularly for larger system sizes.
The dashed line, representing an exponential function $ce^{\alpha n}$ with $c=12$ and $\alpha=0.04$, provides a reference for the growth rate.
For the one-dimensional symmetric and two-dimensional Poisson equations, the increase in iterations is slower than exponential.

\section{Discussion}
\label{sec:discuss}

The error evaluated to calculate the success rate was based on the numerical solution of the SLE.
However, the underlying goal is to solve the original PDE.
The discretization error, which is inherent to numerical methods, plays a crucial role in determining the overall accuracy of the solution. 
In the case of central difference discretization, the error is typically of order $O(h^2)$, where $h=1/(n+1)$ is the mesh size.
To improve the accuracy of the PDE solution, it is necessary to increase $n$, which in turn increases the computational cost.

From an implementation perspective, the system size $n$ is a more critical parameter than the precision parameter $b$.
The number of iterations, and thus, computing time, increases with system size for a given level of precision.
In contrast, the optimal value of $b$ can be determined for a specific system size, as shown in Figs.~\ref{fig:tot-itr} and \ref{fig:itr2}.

Linear-solver algorithms for PDEs typically have computational complexities ranging from $O(n)$ to $O(n^2)$.
In contrast, the proposed method exhibits exponential or subexponential growth in computational time with increasing system size, as illustrated in Fig.~\ref{fig:n_id}.
However, the coefficient $\alpha$ in the exponential function $e^{\alpha n}$ is relatively small, even for the asymmetric Poisson equation. 
This suggests the proposed method, together with using an Ising machine, can outperform conventional methods for certain classes of simple problems.
It is important to note that this analysis focuses on the number of iterations rather than the actual annealing time.
In the simulated annealing method used in this work, the number of sampling steps is proportional to $nb$.
Moreover, solving large-scale problems often requires longer annealing times, as discussed earlier. 

\section{Conclusions}

This study has explored an annealing-based approach for solving PDEs, building upon the algorithm proposed in Ref.~\cite{Krakoff2022} for eigenvalue problems, and analyzed the impact of system size and annealing time on the number of iterations and the success rate.
The number of iterations increases subexponentially or exponentially with system size, with a relatively small exponent.
However, for larger system sizes and shorter annealing times, the success rate can diminish significantly.
The specific problem characteristics also influence computational performance. 
Symmetric problems generally require fewer iterations than asymmetric ones.

Although the presented examples were relatively small-scale problems that can be efficiently solved using traditional methods, the proposed method has the potential for solving larger problems.
By utilizing Ising machines that can handle large-scale problems, the proposed method could offer computational advantages for solving certain classes of PDEs.



 \bibliographystyle{IEEEtran} 
 \bibliography{sle.bib}
\end{document}